\pdfoutput=1\relax
\documentclass[reqno]{amsart}

\usepackage{verbatim}
\usepackage[textsize=scriptsize]{todonotes}
\usepackage{tikz-cd}
\usepackage{etoolbox}
\usepackage{etex}
\usepackage{hyperref}
\usepackage{cleveref}
\usepackage[T1]{fontenc}

\usepackage{chemarr}
\usepackage{amssymb}
\usepackage{amsmath}
\usepackage{comment}
\usepackage{spectralsequences}
\usepackage{cleveref}
\usepackage{mathtools}
\usepackage{rotating}
\usepackage{wrapfig}
\usepackage{outlines}
\usepackage{graphicx}

\theoremstyle{definition}

\newtheorem*{dfn*}{Definition}
\newtheorem*{axm*}{Axiom}
\newtheorem*{ntn*}{Notation}
\newtheorem*{exm*}{Example}
\newtheorem*{exr*}{Exercise}
\newtheorem*{int*}{Intuition}
\newtheorem*{qst*}{Question}
\newtheorem*{rmk*}{Remark}

\theoremstyle{plain}

\newtheorem{thm}{Theorem}

\newtheorem*{thm*}{Theorem}
\newtheorem*{prop*}{Proposition}
\newtheorem*{cor*}{Corollary}
\newtheorem*{lem*}{Lemma}
\newtheorem*{cnj*}{Conjecture}

\let\oldwidetilde\widetilde
\protected\def\widetilde{\oldwidetilde}

\DeclareMathOperator{\Ss}{\mathbb{S}}

\DeclareMathOperator{\F}{\mathbb{F}}

\DeclareMathOperator{\Z}{\mathbb{Z}}

\newcommand{\wt}{\widetilde}

\DeclarePairedDelimiter\abs{\lvert}{\rvert}%
\makeatletter
\let\oldabs\abs
\def\abs{\@ifstar{\oldabs}{\oldabs*}}

\usepackage{tikz}
\usetikzlibrary{matrix,arrows,decorations}
\usepackage{tikz-cd}

\usepackage{adjustbox}

\let\oldtocsection=\tocsection
 
\let\oldtocsubsection=\tocsubsection
 
\let\oldtocsubsubsection=\tocsubsubsection
 
\renewcommand{\tocsection}[2]{\hspace{0em}\oldtocsection{#1}{#2}}
\renewcommand{\tocsubsection}[2]{\hspace{1em}\oldtocsubsection{#1}{#2}}
\renewcommand{\tocsubsubsection}[2]{\hspace{2em}\oldtocsubsubsection{#1}{#2}}

\usepackage{etoolbox}
\newtoggle{draft}
\togglefalse{draft}

\iftoggle{draft} {
\usepackage[margin=1.5in]{geometry}
\newcommand{\NB}[1]{\todo[color=gray!40]{#1}}
\newcommand{\TODO}[1]{\todo[color=red]{#1}}
}{ 
\usepackage[margin=1.5in]{geometry}
\newcommand{\NB}[1]{}
\newcommand{\TODO}[1]{}
\renewcommand{\todo}[1]{}
\renewcommand{\todo}[1]{}
}
\title{An Extension in the Adams Spectral Sequence in Dimension 54}
\date{\today}

\author{Robert Burklund}
\address{Department of Mathematics, MIT, Cambridge, MA, USA}
\email{burklund@mit.edu}

\begin{document}
\begin{abstract}
  We establish a hidden extension in the Adams spectral sequence converging to the stable homotopy groups of spheres at the prime 2 in the 54-stem. 
  This extension is exceptional in that the only proof we know proceeds via Pstragowski's category of synthetic spectra. This was the final unresolved hidden $2$-extension in the Adams spectral sequence through dimension 80. We hope this provides a concise demonstration of the computational leverage provided by $\F_2$-synthetic spectra.
\end{abstract}
\maketitle

The determination of the graded ring of stable homotopy groups of spheres is one of the most concrete and difficult questions in stable homotopy theory. Recent theoretical advances going via motivic homotopy theory beginning with \cite{MR2629898}, \cite{MR3416112} and \cite{GWX} have culminated in \cite{IWX} which constitutes the largest single leap forward in the number of stems we understand, thus far.

Though \cite{IWX} provides near complete information through the 90-stem several uncertainties remain.  The first of these is in the 54-stem where, as pointed out by John Rognes, the argument in \cite{stablestems} regarding the 2-extension contains a mistake which leaves open whether $\kappa \bar{\kappa}^2$ xor $\eta \epsilon \theta_{4.5} + \kappa \bar{\kappa}^2$ is divisible by 2 \cite[Remark 7.11]{IWX}. This article is intended to be the first of several which resolve these uncertainties.

\begin{thm} \label{thm}
  The element $\kappa \bar{\kappa}^2$ in $\pi_{54}$ is divisible by 2.
\end{thm}

In order to prove \cref{thm} we will lift it to a statement in the category of $\F_2$-synthetic spectra where it becomes easier to prove. The category of $\F_2$-synthetic spectra, constructed in \cite{Pstragowski}, is a stable presentably symmetric monoidal $\infty$-category whose objects constitute ``formal Adams spectral sequences'' in the sense of the $E^2$-model-category originally studied by Dwyer, Kan and Stover in \cite{MR1250765}. For a proper introduction to synthetic spectra we direct the reader to \cite{Pstragowski} and to \cite{BHS} for a more computational viewpoint. In this paper we will make use of three main properties of synthetic spectra.

\begin{figure}[t]
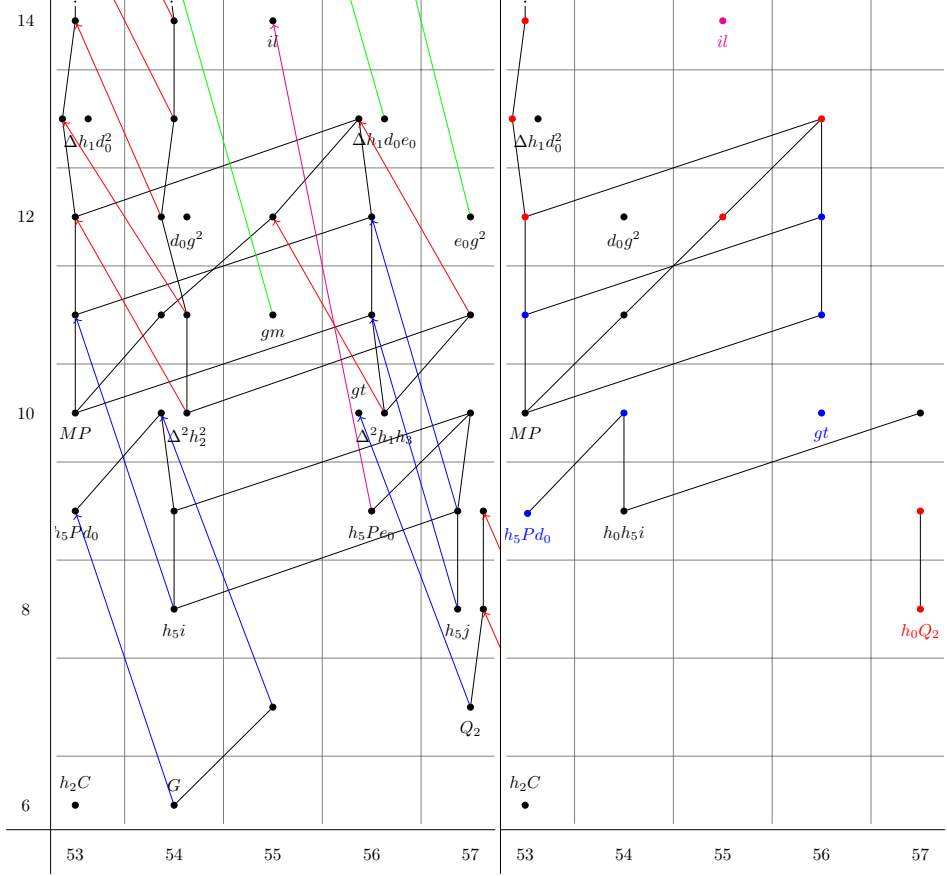

  \centering
  The $\F_2$-Adams spectral sequence and the $\F_2$-synthetic homotopy of the sphere \\
  \scalebox{.65}{
  \printpage[ name = ASS, page = 2 ]
  \printpage[ name = ASS2, page = 2, no y ticks, x axis tail = 0cm ]}
  \caption{ {\footnotesize
      Left: the $\F_2$-Adams spectral sequence near $\kappa \bar{\kappa}^2$.
      Right: the bigraded homotopy groups of the $\F_2$-synthetic sphere in the same region.
      Black dots denote $\tau$-torsion-free classes,
      red dots denote $\tau$-torsion classes,
      blue, green and magenta dots denote $\tau^2$, $\tau^3$ and $\tau^4$ torsion classes respectively.
      In order to reconstruct the group in a given bidegree one must examine all degrees lying above it. For example $\pi_{55,66} \cong \Z/2 \oplus \Z/16$ with generators $\tau^3\{il\}$ and $\tau^{14}\rho_{55}$ respectively. In order to translate from the left picture to the right we remove all boundaries and color the permanent cycles by the length of Adams differential that hits them. This translation is made precise in \cite[Theorem 9.19]{BHS}.} }
  \label{fig1}
\end{figure}

\begin{itemize}
\item Each synthetic spectrum $X$ has bigraded homotopy groups $\pi_{a,b}(X)$ and the monoidal unit $\Ss$ has the property that its bigraded homotopy groups, which we denote simply $\pi_{a,b}$, form a ring. Further, this ring encodes all information present in the Adams spectral sequence in a sense made precise in \cite[Theorem 9.19]{BHS}. A weak form of this correspondence is demonstrated in \cref{fig1}.
\item There is an element $\tau \in \pi_{0,-1}$ with the property that $C\tau$ is a ring and an identification of homotopy rings $\pi_{t-s,t}(C\tau) \cong \mathrm{Ext}_{\mathcal{A}}^{s,t}(\F_2, \F_2)$.
  The unit provides a graded ring map
  $$ \pi_{t-s,t} \to \mathrm{Ext}_{\mathcal{A}}^{s,t}(\F_2, \F_2). $$
\item There is a topological realization to the category of spectra given by inverting $\tau$.
  Again, the unit provides a graded ring map
  $$ \pi_{t-s,t} \to \pi_{t-s}. $$
\end{itemize}


\begin{sseqdata}[ name = ASS, xscale=2.0, yscale=2.0, x range = {53}{57}, y range = {6}{14}, x tick step = 1, y tick step = 2, class labels = {below}, classes = fill, grid = crossword, Adams grading, lax degree]
  \class["h_2C" above](53,6)  
  \class["h_5Pd_0"](53,9)

  \class["MP"](53,10)
  \foreach \i in {11,...,16} { \class(53, \i) \structline }
  \class(54,11) \structline(53,10)
  \class(55,12) \structline
  \class(56,11) \structline(53,10)
  \class(56,12) \structline \structline(53,11)
  \class(56,13) \structline \structline(53,12) \structline(55,12)
  \class["\Delta h_1 d_0^2"](53,13)
  
  \class["G" above](54,6) \class(55,7) \structline

  \class["h_5 i"](54,8) \class(54,9) \structline \class(54,10) \structline \structline(53,9)

  \class["\Delta^2 h_2^2"](54,10)
  \foreach \i in {11,...,15} { \class(54, \i) \structline }

  \class["gm"](55,11)
  \class(54,15)
  \d[green]2(55,11)(54,15) \replacetarget[green]
  \class["d_0g^2"](54,12)
  \class["il"](55,14)
  \class(55,17) \class["\Delta h_1d_0e_0"](56,13) \d[green]2(56,13,2)(55,17)

  \class["h_5 Pe_0"](56,9) \d[magenta]2(56,9)(55,14) \replacetarget[magenta]
  \class(57,10) \structline \structline(54,9)
  \class(57,9) \structline \structline(54,8)
  \class["h_5j"](57,8) \structline 

  \class["Q_2"](57,7) \class(57,8) \structline \class(57,9) \structline
  \class(58,6) \d[red]2(58,6)(57,8,-1) \replacetarget[red]
  \class(58,7) \d[red]2(58,7)(57,9,-1) \replacetarget[red]
  
  \class["gt" above](56,10)
  \class["\Delta^2 h_1h_3"](56,10) \structline(56,11)

  \class(57,11) \structline(57,11)(56,10,2) \structline(57,11)(54,10,2)
  \class["e_0g^2"](57,12) \class(56,16) \d[green]2(57,12)(56,16)

  \d[blue]2(54,6)(53,9) \replacetarget[blue]
  \d[blue]2(55,7)(54,10) \replacetarget[blue]
  \d[blue]2(54,8)(53,11) \replacetarget[blue]
  \d[red]2(54,10,2)(53,12) \replacetarget[red]
  \d[red]2(54,11,2)(53,13) \replacetarget[red]
  \d[red]2(54,12)(53,14) \replacetarget[red]
  \d[red]2(54,13)(53,15) \replacetarget[red]
  \d[red]2(54,14)(53,16) \replacetarget[red]
  
  \d[red]2(56,10,2)(55,12) \replacetarget[red]
  \d[red]2(57,11)(56,13) \replacetarget[red]

  \d[blue]2(57,7)(56,10) \replacetarget[blue]
  \d[blue]2(57,8)(56,11) \replacetarget[blue]
  \d[blue]2(57,9)(56,12) \replacetarget[blue]
  
\end{sseqdata}

\begin{sseqdata}[ name = ASS2, xscale=2.0, yscale=2.0, x range = {53}{57}, y range = {6}{14}, x tick step = 1, y tick step = 2, class labels = {below}, classes = fill, grid = crossword, Adams grading, lax degree]
  \class["h_2C" above](53,6)  
  \class[blue, "h_5Pd_0", below right](53,9)

  \class["MP"](53,10)
  \class[blue](53,11) \structline
  \foreach \i in {12,...,16} { \class[red](53, \i) \structline }
  \class(54,11) \structline(53,10)
  \class[red](55,12) \structline
  \class[blue](56,11) \structline(53,10)
  \class[blue](56,12) \structline \structline(53,11)
  \class[red](56,13) \structline \structline(53,12) \structline(55,12)
  
  \class["h_0h_5i"](54,9) \class[blue](54,10) \structline \structline(53,9)
  \class["\Delta h_1d_0^2"](53,13)

  \class["d_0g^2"](54,12)
  \class[magenta, "il"](55,14)
  
  \class(57,10) \structline(54,9)
  
  \class[red, "h_0Q_2"](57,8) \class[red](57,9) \structline  
  \class[blue, "gt"](56,10)
  
\end{sseqdata}

The maps we have described provide a more direct link between topoology and ext over the Steenrod algebra than previous techniques allowed. The main reason we use $\F_2$-synthetic spectra in this paper over other choices of Adams type homology theory is the availability of large-scale machine calculations of the ext ring of the Steenrod algebra from \cite{E2}.

Breaking with the notation from \cite[Theorem 9.19]{BHS} we will use the following notation for elements of $\F_2$-synthetic homotopy:
Given a permanent cycle $a$, what previously would have been called $\wt{a}$ we now call $\{a\}$.
Given a element $\alpha \in \pi_*$, what previously would have been called $\wt{\alpha}$ we now call $\alpha$, with the exception of $\wt{2} \in \pi_{0,1}$.
These changes were suggested by Isaksen and Xu.
We hope they better match with existing conventions in the literature.

\begin{proof}
  We will prove the theorem by showing that $\wt{2} \tau^2 \{h_0h_5i\} = \kappa \bar{\kappa}^2 \tau^4$.
  Using \cite[Theorem 9.19(4)]{BHS} and the Adams spectral sequence calculations from \cite{IWX} we learn that there are unique coefficients $a_i \in \{0,1\}$ such that,
  \begin{align*}
    \wt{2} \{h_0h_5i\}
    &= a_1 \eta \{h_5Pd_0\} + a_2 \eta \{MP\} \tau + a_3 \kappa \bar{\kappa}^2 \tau^2 \\
    \bar{\kappa} \{h_0h_5i\}
    &= a_4 \nu \{\Delta^2 h_2g\} \tau + a_5 \eta \{MP\} \bar{\kappa} \tau^2 + a_6 \kappa \bar{\kappa}^3 \tau^3 \\
    \wt{2} \bar{\kappa} \{h_0h_5i\}
    &= a_7 \nu \{\Delta^2 h_2g\} + a_8 \eta \{MP\} \bar{\kappa} \tau + a_9 \kappa \bar{\kappa}^3 \tau^2
  \end{align*}

  \begin{itemize}
  \item Since $h_5Pd_0$ is hit by a $d_3$ differential, $\{h_5Pd_0\}\tau^2 = 0$.
  \item In $C\tau$ we can read off from \cite{E2} that $h_0^2h_5i = h_1Pd_0h_5$. This implies $a_1 = 1$.
  \item From \cite[Lemma 7.21]{IWX} we know that $2\{h_0h_5i\}$ is nonzero after inverting $\tau$.
    Thus, at least one of $a_2$ and $a_3$ is nonzero.
  \item From \cite[Theorem 9.19]{BHS} and \cite{IWX} we can compute that $\pi_{73,86}$ has no $\tau$-torsion, therefore $\bar{\kappa} \{h_5Pd_0\} = 0$.
    This allows us to conclude that $a_7 = 0$, $a_2 = a_8$ and $a_3 = a_9$.
  \item Similarly, we can compute that $\pi_{74,90} \cong \Z/2 \oplus \Z/2$, therefore $2 \kappa \bar{\kappa}^3 = 0$.
  \item Using the fact that one of $a_8$ and $a_9$ is nonzero we learn that
    $$ 0 \neq \wt{2} \bar{\kappa} \{h_0h_5i\} = \wt{2} (a_4 \nu \{\Delta^2 h_2g\} \tau + a_5 \eta \{MP\} \bar{\kappa} \tau^2 + a_6 \kappa \bar{\kappa}^3 \tau^3) = a_4 \wt{2} \nu \{\Delta^2 h_2g\} \tau, $$
    which implies $a_4 = 1$.
  \item In $C\tau$ we can read off from \cite{E2} that $h_0h_2 \Delta^2h_2g = 0$. This implies $a_8 = 0$. \footnote{We follow \cite[Table 4]{IWX} in using the notation $\Delta^2h_2g$ for Bruner's element $x_{13,34}$.}
  \end{itemize}

  Altogether, we may conclude that,
  \begin{align*}
    \wt{2}\{h_0h_5i\} &= \eta \{h_5Pd_0\} + \kappa \bar{\kappa}^2 \tau^2 \qedhere
  \end{align*}
  
\end{proof}

We end by commenting on why this argument cannot be run in the ordinary category of spectra.
The key point is that $\eta \epsilon \theta_{4.5} \bar{\kappa} = 0$ which means that after multiplying by $\bar{\kappa}$ we've lost the ability to distinguish between the two possible extensions. In the Adams spectral sequence $h_1g MP$ is hit by a $d_4$ differential which translates into the product being $\tau^3$-torsion, but still non-trivial synthetically. Then, we are able to resolve the extension because its source jumps filtration in a way that makes it possible to read $a_8$ off from homological algebra. This same strategy is not viable in the $\mathbb{C}$-motivic world because the differential killing $Mh_1 d_0^2$ is shorter than the extension.

\subsection*{Acknowledgments}
We thank Dan Isaksen for comments on a draft.
We thank Zhouli Xu for comments, suggesting this question and for introducing us to the joy of computational stable homotopy theory. 

\bibliographystyle{alpha}
\bibliography{bibliography}

\end{document}